\documentclass[12pt]{amsart}

\usepackage{PTDefn}

\title{Representation Stability for the Kontsevich space of stable maps}

\date{\today}
\author{Philip Tosteson}
\address{Department of Mathematics, University of Chicago, Chicago, IL}
\email{\href{mailto:ptoste@umich.edu}{ptoste@math.uchicago.edu}}
\urladdr{\url{https://math.uchicago.edu/~ptoste/}}
\thanks{This work was supported by NSF grant  DMS-1903040.}

\newcommand{\mainbound}{(13i/2 + 1)(i + 4 g + 4 \deg \alpha + 1)}

\begin{document}
\maketitle
\begin{abstract}
	For an algebraic variety $X$, curve class $\alpha \in \rN_1(X)$, and genus $g \in \bbN$ we consider the sequence of $\rS_n$ representations obtained from the homology of the Kontsevich space of stable maps to $X$,  $\overline \cM_{g,n}(X,\alpha)$.   Using the category of finite sets and surjections, we prove a representation stability theorem that governs the behavior of this sequence of representations for $n$ sufficiently large.  
\end{abstract}

\section{Introduction}
Let $X$ be a projective algebraic variety over $\bbC$, and let $\alpha \in \rN_1(X)$ be a numerical equivalence class that is a sum of curve classes.   Let  $\bMgn(X,\alpha)$  be the Kontsevich space of stable maps to $X$.  The moduli space $\bMgn(X,\alpha)$  parameterizes genus $g$  nodal curves $C$ equipped with $n$ marked points and a map $f: C \to X$,  satisfying (1) that $f_*[C] = \alpha$ and (2) any component of $C$ contracted by $f$ has at least $3$ marked or nodal points.   

In this paper, we fix $i,g$ and $\alpha$  and consider the asymptotic behavior of the rational homology  of $\bMgn(X,\alpha)$  for  $n \gg 0$.    The following theorem is an application of our main result.

\begin{thm}\label{thm:consequences}
	Let $X$ be a projective variety over $\bbC$.  Fix $i,g \in \bbN$ and $\alpha \in \rN_1(X)$, and let $C = (13i/2 + 1)(i + 4 g + 4 \deg \alpha + 1)$. \begin{enumerate} 		 \item The generating function $$\sum_{n \in \bbN} \dim H_i(\bMgn(X,\alpha), \bbQ) t^n$$  is rational with denominator  $\prod_{j = 1}^C (1- jt)^{e_j}$ for some $e_j \in \bbN$.    In particular, $\dim H_i(\bMgn(X,\alpha)$  equals a sum of polynomials times  integer exponentials for $n \gg 0$. 
\\ 
	\item  For $n \gg 0$, the function $\sum_{n \in \bbN} \dim H_i(\bMgn(X,\alpha), \bbQ)^{\bS_n} t^n$  is rational.  The poles of the denominator are at roots of unity of order less than $C$.    
\end{enumerate}
\end{thm}

Theorem \ref{thm:consequences} follows from a structural statement about the homology of $\bMgn(X,\alpha)$, which we now describe.   

Let $\FS$ be the category whose objects are the finite sets $[n] = \{1,\dots, n\}$  for $n \in \bbN$, and whose morphisms are surjections.  An $\FSop$ module is a contravariant functor  from $\FS$ to the category of $\bbQ$ vector spaces.

\begin{defn}
Let $\bP_d$ be the projective $\FSop$ module $[n] \mapsto \bbQ \FS([n],[d])$.  
  We define the subcategory of $\FSop$ modules of \emph{height $\leq C$} to be the smallest subcategory containing  $\bP_d$ for $d \leq C$ that is closed under isomorphism,  submodules, quotient modules, and extensions.  
\end{defn}
 
\begin{rmk} The notion of a height $\leq C$ module appears  in \cite{PR}.   There  $\FSop$ modules of height $\leq C$ are called  $C$-smallish.   \end{rmk}

\begin{thm}\label{mainthm}
   Let $\alpha \in \rN_1(X)$, $i,g \in \bbN$.  Then there is an $\FSop$ module whose underlying sequence of $\rS_n$ representations is $n \mapsto H_i(\bMgn(X,\alpha))$.  This $\FSop$ module has height $\leq \mainbound$.   
\end{thm}

Theorem \ref{mainthm}  together with structural results on $\FSop$ modules due to Sam--Snowden \cite{sam2017grobner} and Tosteson \cite{tosteson2021categorifications}  immediately imply Theorem \ref{thm:consequences}.   More broadly, the structure theory for $\FSop$ module characters developed in \cite{tosteson2021categorifications} applies to the sequence of $\rS_n$ representations $n \mapsto H_i(\overline \cM_{g,n}(X,\alpha))$.  In particular, the Frobenius character of this sequence of representations  is a rational function.

The strategy of the proof of Theorem \ref{mainthm} is  similar to the case where $X$ is a point, established in \cite{tosteson2018stability}, but we encounter several additional difficulties.  Combinatorially, the stratification of by graphs $\overline \cM_{g,n}(X,\alpha)$ is more complicated because of the labeling of vertices by elements of $\rN_1(X)$.  Geometrically, the argument in \cite{tosteson2018stability} used that the strata of $\overline\cM_{g,n}$  are products of  moduli spaces $\cM_{h,m}$.  This is no longer true for Kontsevich spaces: to handle this failure we introduce the notion of a \emph{free half-edge} and we prove new combinatorial statements about free half edges.    

Most significantly, the proof of \cite{tosteson2018stability} used a purity argument,  which is not available for Kontsevich spaces because they are not smooth.  In \S \ref{genus0} we resolve the locus  $\overline \cM_{g,n,0} \subseteq \bMgn$  of stable curves whose irreducible components all have genus zero.  Our resolution yields a chain complex that computes $H_*(\overline \cM_{g,n,0})$ in terms of the homology of \emph{matroid independence complexes} associated to genus $g$ graphs (and $H_*(\bMon)$).   This complex is the main novelty of our paper, and may be of independent interest.    

Finally, we note that even in the in the case where $X$ is a point, the bound of Theorem \ref{mainthm}  is stronger than the statement proved in \cite{tosteson2018stability}, due to our sharpening of several combinatorial arguments  in the process of generalizing them.

\subsection{Acknowledgements} I would like to thank Ming Zhang, for suggesting I consider representation stability for Kontsevich mapping spaces.    I would also like to thank Dan Petersen for pointing out that an effective class $\alpha \in \rN_1(X)$ can only be written as a sum of effective classes in finitely many ways, and suggesting that it would be possible to extend the methods of \cite{tosteson2018stability} to Kontsevich mapping spaces.

\subsection{Notation and Conventions}  Throughout the body of this paper all homology and cohomology groups are implicitly taken with $\bbQ$ coefficients.

		We work over $\bbC$ and   $\overline \cM_{g,n}(X,\alpha)$  denotes the coarse moduli space associated to the Kontsevich mapping stack, as defined by Behrend--Manin \cite{behrend1996stacks} (also see \cite{fulton1996notes}).   Since the rational (co)homology groups of a Deligne--Mumford stack are isomorphic to those of its associated coarse moduli space, the homology of the stack also carries $\FSop$ module structure and Theorem \ref{mainthm} holds for it.  All of our geometric arguments have analogs for the moduli stacks, which can be made using the functor of points.  A reader who prefers to work with stacks may translate accordingly.

\section{Background on $\FSop$ modules and definition of action}
	The category of $\FSop$ modules (as defined in the introduction) is canonically equivalent to the category of contraviariant functors from the category of \emph{all} finite sets and surjections to $\Vec_\bbQ$.   We will pass between these two notions interchangeably.  

We recall the category $\BT$ of binary trees (see \cite[Definition 2.1]{tosteson2018stability}).   The set of objects of $\BT$ is $\{[n]\}_{n \in \bbN}$.  The set of morphisms from $[m]$ to  $[n]$,  denoted $\BT([m],[n])$, is the set of binary trees with roots labelled by $[n]$ and leaves labeled by $[m]$.   
There is an action of $\BT\op$  on $\bMgn(X,\alpha), n \in \bbN$;  More precisely, there is a functor from $\BT\op$ to algebraic varieties, given on objects by $[n] \mapsto \bMgn(X,\alpha)$,  with morphisms defined as follows.  

We recall from   \cite[Definition 2.5]{tosteson2018stability} that for every binary forest  $F \in \BT([m],[n])$  there is canonical variety $L_F$  marked by $[m] \sqcup [n]$,  which is a disjoint union of stable genus $0$ curves with dual graph $F$.  Then  $F$ determines a morphism $\bMgn(X,\alpha) \to \overline \cM_{g,m}(X,\alpha)$  by $C \mapsto C  \sqcup_{[n]} L_F$.

\begin{prop}
	The action of $\BT\op$ on $\bMgn(X,\alpha)$ induces an action of $\FSop$ on $H_i(\bMgn(X,\alpha))$ for every $i \in \bbN$.  
\end{prop}
\begin{proof}
	The proof of \cite[Proposition 2.7]{tosteson2018stability} directly extends to this case.  Namely,  any two forests  $F_1, F_2 \in \BT([m],[n])$ which have the same image $g \in \FS([m],[n])$  correspond to points in $\prod_{i \in [n]} \overline \cM_{0,g\inv(i)}$.  Since these points are connected by a path,  the two forests induce the same map on homology.  
\end{proof}

\subsection{Closure properties of height}

We will need the following elementary closure properties of $\FSop$ modules with height $\leq i$.  They were originally proved for $i$-small $\FSop$ modules in  \cite{proudfoot2017configuration}, and they easily extend  to $\FSop$ modules  of height $\leq i$ (also called $i$-smallish $\FSop$ modules).  

\begin{prop}\label{shiftclosure}
  Let $M$ be an $\FSop$ module of height $\leq i$.  Let $k \in \bbN$.  Then the $\FSop$ module $\Sigma^{k}  M$ defined by $X \mapsto  M_{[k] \sqcup X}$ has height $\leq i$.  
\end{prop}
\begin{proof}
	We may reduce to the case  $M = \rP(i')$ for $i' \leq i$.  Then see Lemma 4.3 of \cite{proudfoot2017configuration}.
\end{proof}

\begin{prop}\label{convolutionclosure}
	Let $N^1, \dots, N^r$  be $\FSop $ modules which are respectively of height less than or equal to $i_1, \dots, i_r$.  Then the $\FSop$ module $\oast_{j = 1}^r N^j$  defined by $$X \mapsto \bigoplus_{f: X \to [r]}\bigotimes_{j = 1}^r N^j_{f \inv(j)},$$ has height $\leq \sum_{j = 1}^r i_j$.  
\end{prop}
\begin{proof}
	We may reduce to the case $N^j = \rP(i_r')$ for $i_r' \leq i_r$.  Then  see Lemma 4.2 of \cite{proudfoot2017configuration}.
\end{proof}

	By definition modules of height $\leq h$ are closed under kernels, cokernels and extensions  so the following proposition is immediate. 
\begin{prop}\label{SSclosure}
 	Let $E^1_{p,q} \implies M^{p+q}$ be a convergent spectral sequence of $\FSop$ modules.    If there are finitely many $p,q \in \bbZ$ with $p + q = j$  and for every such $p,q$ the module $E^1_{p,q}$ has height $\leq i$   then $M^j$ has height $\leq i$.  
\end{prop}

\section{Stratification of $\overline \cM_{g,n}(X,\alpha)$}  

The Kontsevich mapping space has a standard stratification by graphs equipped with certain auxiliary data.   Here we define a coarsening of this stratification, such that $\FSop$ acts on the homology groups of the associated spectral sequence.  

\subsection{Posets of stable decorated graphs}

\begin{defn}
	A \emph{graph} $G$  consists of a finite set $H_G$ of half-edges, an involution $\sigma_G: H_G \to H_G$  and a partition of $H_G$.  
The blocks of the partition is the set of \emph{vertices} of $H$,  denoted $V_G$.  
 An \emph{internal edge} of $G$ consists of an orbit of $\sigma_G$ of size two.  
A half-edge fixed by $\sigma_G$ is called an \emph{external edge}. 
 We say that $v \in V_G$ is \emph{adjacent} to a half-edge $e \in H_G$ if $e \in v$.     
We write $n(v)$ for the number of half-edges adjacent to $v$, in other words  the \emph{valence of $v$}.  
 We say that the vertices $v,w$ are adjacent if there exist half-edges  $e \in v, f \in w$ such that $\sigma_G(f) = e$.  
\end{defn}

\begin{remark}
We will sometimes abuse notation and write $n(v)$ to denote the \emph{set} of half edges adjacent to $v$.     In particular we will use,  $\cM_{g,n(v)}$  to denote the moduli space of genus $g$ curves  with marked points labelled by the half--edges adjacent to $v$. (More formally,  the set of edges adjacent to $v$ is simply $v$ so we should write $\cM_{g,v}$,  but this notation may cause confusion).
\end{remark}

An \emph{isomorphism} of graphs $f: G \to G'$  is a bijection $H_G \to H_{G'}$  preserving the partition and involution. 
 So we may speak of the automorphism group of a graph, and the set of isomorphism classes of graphs.  

A \emph{contraction} of graphs $\pi: G \to G'$  is an embedding $\iota : H_{G'} \to H_{G}$ preserving the involution, which is a bijection on fixed points and satisfying the property 
\begin{itemize} \item for every edge $\{h_1, h_2\} \subseteq H_{G}$ in the complement of $\im \iota$, if  $j_1, j_2 \in H_{G'}$  are such that $\iota(j_i)$ and $h_i$ are adjacent to the same vertex  for $i = 1,2$,  then  $j_1$ and $j_2$ are adjacent to the same vertex.  
\end{itemize}
Topologically, the map $\pi: G \to G'$  corresponds to the map that contracts the set of edges in the complement of the image of $\iota$ and permutes the remaining edges according to the bijection induced by $\iota$.  In particular, the contraction $\pi$ induces a morphism $ V_G \to V_{G'}$, which we will also denote by $\pi$.

\begin{defn}
Let $h \in \bbN$ and $\beta \in \rN _1(X)$.   A \emph{stable decorated graph of genus $h$ and class $\beta$}  consists of a connected graph  $G$ and functions $g:  V_G \to \bbN$  and $\alpha: V_G \to \rN_1(X)$,   subject to the following conditions:
\begin{itemize}
\item  $\sum_{v \in V_G} \alpha(v) = \beta$ and $h^1(G) + \sum_{v \in V_G} g(v) = h$,   where $h^1(G)$ is the first betti number of $G$.  
\item   If $\alpha(v) = 0$ and $g(v) = 0$ then $n(v) \geq 3$
\item If $\alpha(v) = 0$ and $g(v) =1$ then $n(v) \geq 1$.  
\item  For every $v \in V_G$ the class $\alpha(v)$ is a sum of curve classes $[C] \in \rN_1(X)$.
\end{itemize}
We write $n(G)$ for the number of external edges of $G$.   We say that a vertex $v \in G$ is \emph{undecorated} if $\alpha(v) = 0$,  and that it is \emph{plain} if $g(v) = 0$ and $\alpha(v) = 0$.
\end{defn}

An isomorphism of stable decorated graphs is an isomorphism of graphs that preserves $g$ and $\alpha$.    A contraction of stable graphs $\pi: G \to G'$ is a contraction of graphs satsifying  $\alpha(v') = \sum_{v \in \pi \inv(v')}  \alpha(v)$ and $g(v') =\sum_{v \in \pi \inv(v')} g(v)$.   

\begin{defn}
We write $\Stab(h,n,\beta)$ for the set of isomorphism classes of stable decorated genus $h$ graphs with $n$ external edges, and class $\beta$. Then $\Stab(h,n,\beta)$ is naturally a poset,  where we define $[G] \leq [G']$ if there is a contraction $G \to G'$. 
\end{defn}

\begin{prop}\label{prop:finite}
	The set $\Stab(h,n,\beta)$ is finite.  
\end{prop}
\begin{proof}
 	 Let $L$ be an ample divisor on $X$.   Let $G \in \Stab(h,n,\beta)$.   There are at most $L \cdot \beta + h + n$ many vertices in $G$ that are adjacent to an external edge or not plain,  and they contribute at most $L \cdot \beta + h + n$ to $\chi(G)$.  A plain vertex not adjacent to an external edge is at least trivalent, so it contributes $ \leq -1/2$ to  $\chi(G) \geq 1 - h$.  This shows that the number of vertices of $G$ is uniformly bounded, hence  there are finitely many isomorphism classes of graph   in $\Stab(h,n\beta)$ (neglecting the labellings $g(v)$ and $\alpha(v)$).  

The set of curve classes $\alpha$ on $X$ such that $L \cdot \alpha \leq L \cdot \beta$ is finite (see \cite[Example 1.4.31]{lazarsfeld}).   So there are at most finitely many ways to write $\beta$ as a sum of curve classes on $X$.   Hence there are at most fintiely many choices for $\alpha(v)$ and $g(v)$ for a given isomorphism class of graph. 
\end{proof}

\begin{remark} Alternatively, Proposition \ref{prop:finite} follows from the fact that Kontsevich mapping space $\overline \cM_{g,n}(X, \beta)$ is projective and admits a stratification with strata labelled by $\Stab(g,n,\beta)$ (see \S \ref{subsec:strat}).  
\end{remark}

\begin{defn}
	We say that a stable decorated graph $G \in \Stab(h,n,\beta)$ is \emph{saturated}  if there are no adjacent  plain vertices  of $G$.   We write $\rQ(h,n,\beta)$ for the subset of $\Stab(h,n,\beta)$ consisting of saturated graphs.  
\end{defn}

Given a graph $G \in \Stab(h,n,\beta)$, we may associate to it a saturated quotient graph, as follows.

\begin{defn}
  For $G \in \Stab(h,n,\beta)$, choose a spanning forest for the induced subgraph on the set of plain vertices of $G$.   We define $\overline G$ be graph obtained by quotienting the edges of this forest.     The graph $\overline G$ is saturated, and well defined up to isomorphism  (independent of the choice of spanning forest).   
\end{defn}

\begin{prop}\label{prop:poset}
There is a poset structure on $\rQ(h,n,\beta)$ such that set map $$\Stab(h,n,\beta) \to \rQ(h,n,\beta), ~G \mapsto \overline G$$ is a poset morphism.  
\end{prop}
\begin{proof}
We define the relation $\leq_Q$ on $\rQ(h,n,\beta)$ as follows.   We say that $G \leq_Q H$  if there is a sequence of graphs  $G = G_0, G_1, \dots, G_n = H$  with $G_i \in \Stab(h,n, \beta)$  such that  for all $i = 0, \dots, n-1$ one of the following holds:
\begin{enumerate}
\item[(I)]  $G_i \prec G_{i+1}$ 
\item[(II)] $G_i \succ G_{i+1}$ and $G_i$ is obtained from $G_{i+1}$ by collapsing an edge between two distinct plain vertices.
\end{enumerate}
We say that $G_i$ and $G_{i+1}$ are related by a \emph{forward} collapse in case (I)  and by a \emph{backward} collapse in case (II).  
By construction the relation $\leq_Q$ is transitive and reflexive.  We need to show that it is anti-symmetric.    

Suppose that $G \leq_Q H$ and $H \leq_Q G$.  Then there exists a sequence  $G = G_0, \dots, G_s = H,  G_{s+1}, \dots, G_{t} = G$  of graphs related by forward and backwards edge collapses.

We use an invariant $I: \Stab(h,n,\beta) \to \bbZ^3$.  Define $$I(J) = \left(\sum_{v \in J} g(v), ~ - \#\{v \in J \text{ not plain}\}, ~ \sum_{v \in H, \text{ not plain }} n(v)\right).$$  We totally order $\bbZ^3$ lexicographically: define  $(n_1, n_2, n_3) < (n_1', n_2', n_3')$ if for the smallest $i \in \{1,2,3\}$ such that $n_i \neq n_i'$  we have $n_i < n_i'$.  


Suppose that $J \in \Stab(h,n, \beta)$,  $e$ is an edge of $J$,  and $J'$ is quotient of $J$ obtained by collapsing $e$.  We determine the relationship between $J$ and $J'$  in the following (exhaustive) list of cases.

\begin{enumerate}
\item  $e$ is an edge between distinct plain vertices.  Then $I(J) = I(J')$
\item $e$ is a self edge.  In this case, the total  genus increases from $J$ to $J'$, so $I(J)< I(J')$.
\item $e$ is an edge from a plain vertex to a non-plain vertex.   Then $I(J)$ and $I(J')$ only differ in the last entry.   Since the plain vertex has valence $m \geq 3$, we have that $I(J')_3 = I(J)_3 + m - 2   > I(J)_3.$ So $I(J) < I(J')$.
\item  $e$ is an edge between two distinct non-plain vertices.   Then the total genus of $J$ and $J'$ is the same,  and the number of non-plain vertices decreases,  so $I(J) < I(J')$.
\end{enumerate}

From these four cases, it follows that $I$ is a poset morphism.  By this, together with case $(1)$, we have that $I(G_i) \leq I(G_{i+1})$ for all $i$. Since $I(G) = I(G_0) = I(G_t)$, we must have $I(G_i) = I(G_{i+1})$ for all $i$.  

Thus by the case analysis above, if $G_i$ and $G_{i+1}$ are related by a forward edge collapse then that edge collapse must be between distinct genus $0$ vertices.  Therefore $\overline G_i = \overline G_{i+1}$  and so $\overline G = \overline H$.  Since $H$ and $G$ are already saturated we have $G = H$, proving antisymmetry.       

Finally the map $J \mapsto \overline J$ is a poset morphism $\Stab(h,n, \beta) \to \rQ(h,n,\beta)$ by the definition of $\leq_\rQ$.
\end{proof}

\subsection{Stratification}\label{subsec:strat}
  We recall that there is a stratification of $\bMgn(X,\beta)$ by the poset of stable graphs $\Stab(g,n,\beta)$, defined in terms of the dual graph.   

\begin{defn} The \emph{dual graph} of a stable map  $(C,f) \in \bMgn(X,\beta)(\bbC)$  is the graph $G_C \in \Stab(g,n,\beta)$ with \begin{itemize}
\item  a half edge for every marked point of $C$,  which is fixed by $\sigma_{G_C}$
\item  two half edges for every nodal point $p$ of $C$ that are interchanged by $\sigma_{G_C}$, corresponding to the two pairs $(p,C_i)$  where $C_i$ is an irreducible component of $C$ containing $p$
\item a vertex for each irreducible component $C_0$ of $C$,  which is adjacent to every half edge corresponding to a marked point of $C_0$, and to every half edge associated to a pair $(p,C_0)$.
\end{itemize}
The function $g: V_{G_C} \to \bbN$ takes a vertex $v$  to $g(C_v)$,  the genus of its corresponding irreducible component.  Similarly we define $\alpha(v)$  to be  the pushforward of the fundamental class of $C_v$ along $f|_{C_v}$.  
\end{defn}

\begin{defn}
	Given $H \in \Stab(g,n,\beta)$  there is a  closed subscheme $\overline \cM_{H} \subseteq \overline\cM_{g,n}(X,\beta)$,  parameterizing stable maps  maps whose dual graph is $\leq H$ in the partial order on $\Stab(g,n,\beta)$.    Then  $\cM_{H} := \overline \cM_{H} - \bigcup_{K < H} \overline \cM_{K}$  is a locally closed subscheme parameterizing stable maps whose dual graph is $H$.    
\end{defn}

Using the morphism $\Stab(g,n, \beta) \to \rQ(g,n,\beta)$, we may define a stratification of $\bMgn(X,\beta)$ by $\rQ(g,n, \beta)$  by defining  $$\rS(G) := \bigcup_{H \in \Stab(g,n,\beta),\overline H = G}  \cM_H,$$ where $\cM_H$ is the stratum of $\bMgn(X, \beta)$ associated to the stable decorated graph $H$.    

\begin{defn}Given $G \in \rQ(g,n,\beta)$ and $v \in G$,   we define $\cM(v)$ casewise  as follows.    
\begin{enumerate}
\item\textbf{If $v$ is plain:}     Let $e(v)$ number of self edges of $v$  and let $r(v) \subseteq n(v)$ be the set of half-edges are not part of a self edge.   We define $$\cM(v) \subseteq \overline \cM_{e(v),r(v)}(X,0) =  \cM_{e(v),r(v)} \times X$$ to be the moduli space parameterizing stable genus $e$ curves with each  irreducible component is genus $0$, together with a constant map to $X$. \\
\item \textbf{If $v$ is not plain:}  We define $\cM(v)$ to be $\cM_{n(v),g(v)}(X, \alpha(v))$.  
\end{enumerate} 
\end{defn}

\begin{defn}
Let $$\prod_{v \in V_G} \cM(v)  \to  X^{\tilde H_G}$$ be the canonical evaluation map  where $$\tilde H_G := \{h \in H_G ,~  h \text{ is not part of an external or self edge}\}.$$  Define $\tilde \rS(G)$ to be the fiber product of this map with canonical embedding   $X^{\tilde E_G} \to X^{\tilde H_G}$, where $$\tilde E_G:= \{ e \in E_G,~ e \text{ not a self edge of a plain vertex}\}.$$  The space $\tilde S(G)$ parameterizes stable algebraic maps $f_v:  C_v \to X$ for $v \in V_G$, where $f_v$ has class $\alpha(v)$ and the following hold:
\begin{itemize}
\item  when $v$ is  not plain, $C_v$ is smooth of genus $g(v)$ and has marked points   labeled by $n(v)$,  the set of half edges of $v$;
\item  when $v$ is plain, $C_v$ has genus $e(v)$ and is has marked points labeled by $r(v)$,  the set of half edges of $v$ that are not part of self edge;

\item for every edge $e \in E_G$ which is not a self edge of a plain vertex,  the two marked points $p_{1} \in C_{v_1}$ and $p_{2} \in C_{v_2}$  corresponding to half edges of $e$  satisfy $f_{v_1}(p_{v_1}) = f_{v_2}(p_{v_2})$.

\end{itemize}
\end{defn}

The automorphism group of $G$ acts on $\tilde S(G)$.  An element $\sigma \in \Aut(G)$  acts by taking $(C_v)_{v \in G}$  to $(C_\sigma(v))_{v \in G}$  and relabeling the marked points of $C_v$ according to the action of $\sigma$ on half edges.

\begin{prop}\label{quotienthomeo}
	There is a homeomorphism of moduli spaces $\tilde S(G)/\Aut(G) \to S(G)$.  
\end{prop}
\begin{proof}
	We show that $\tilde S(G) \to S(G)$ factors through $\tilde S(G)/\Aut(G)$  and that the induced map is bijective and proper.  

 For $\sigma \in \Aut(G)$, there is a canonical isomorphism between ${\rm glue}(C_v)_{v \in V_G}$ and ${\rm glue} (C_{\sigma(v)})_{v \in V_G}$  so the morphism on coarse moduli factors through  $S(G)/\Aut(G)$.

	The morphism is bijective since if $\glue((C_v)_{v \in V_G})$ is isomorphic to $\glue((C'_v)_{v \in   G})$,   then that isomorphism yields a relabelling $\tau \in \Aut(G)$ of the components and marked points of $(C_v)_{v \in V}$ such that $C_{\tau(v)}  \iso C'_{v}$.  

 To show properness we consider the compactified gluing map $$gl: \prod_{v \in V_G} \overline \cM(v) \times_{  X^{\tilde H_G}} X^{\tilde E_G} \to \overline {S(G)},$$  where we define  $\overline \cM(v)$ to be $\cM(v)$ if $v$ is plain  and $\overline \cM(v) = \overline \cM_{g(v),n(v)}(X, \alpha(v))$  otherwise.  

 We claim $f\inv( S(G)) = \tilde S(G).$   Since by definition $f\inv(S(G)) \supseteq \tilde S(G)$,  it suffices to show that for any curve $(C_v)_{v \in V_G}$ in the complement of $\tilde S(G)$  that $f((C_v)) \not \in S(G)$.  Let $J$ be the dual graph of $f((C_v))$.  For at least one non-plain vertex vertex $v_0 \in V_G$,  we have that $C_{v_0}$ lies in the boundary of $\overline \cM(v)$.   Therefore $G$ is obtained from $J$ by contracting edges  with the final edge contraction (corresponding to $v_0$) falling into one of the cases (2),(3),(4)  listed in the proof of Proposition \ref{prop:poset}.   Thus $I(J) < I(G)$  and therefore $f(C) \not \in S(G)$  (because if we had $\overline J = G$  then we would have $I(J) = I(G)$).  

Since $f$ is proper, its restriction to $f \inv(S(G)) = \tilde S(G)$ is as well.   Hence the map $\tilde S(G) \to \tilde S(G)/\Aut(G) \to S(G)$ is proper.   
\end{proof}


\subsection{Free half-edges}

\begin{defn}
Let $G \in \rQ(g,n)$.  We say that a half-edge $f \in H_G$ is \emph{free} if the vertex $v$ it is adjacent to is not plain, and one of the following holds:
\begin{itemize}
\item $v$ is undecorated
\item $v$ is decorated and every path from $f$ to a half edge that is adjacent to a decorated vertex passes through $v$.
\end{itemize}  
For $v \in G$, we write $F(v)$ for the number of free half-edges adjacent to $v$.  We say that a half-edge $f$ is \emph{bound} if it is not adjacent to a plain vertex and not free.
\end{defn}

Our main motivation for this definition is the following proposition.  

\begin{prop}\label{freebound}
	Let $\alpha \in \rN_1(X)$ and let $L \in \rN^1(X)$ be ample. Let $G \in \rQ(g,n,\alpha)$.  Suppose that the number of free half-edges of $G$ is $> i + g + L \cdot \alpha$.  Then $H_i^{\rm BM}\left(\tilde\rS(G) \right) = 0$.  
\end{prop}

Before proving Proposition \ref{freebound}, we first describe $\tilde \rS(G)$  using free edges. 

\begin{defn}
 For $v \in G$, we define $\cM'(v)$ as follows. 
\begin{itemize}
\item   If $v$ is decorated, we let $\cM'(v) = \cM(v)$.

\item
If $v$ is not decorated and $g(v) \geq 1$,  we let $\cM'(v) = \cM_{n(v),g(v)}$.
\item  If $v$ is plain, then $\cM'(v) \subseteq  \overline \cM_{e(v),r(v)} $  is the locus parameterizing stable genus $e$ curves with points marked by $r$ such that every irreducible compontent is genus $0$.
    
\end{itemize}
\end{defn}

Next we define an equivalence relation on the bound half-edges of $G$,  by declaring that $h_1 \sim h_2$ if $h_1$ and $h_2$ are connected by a path in $G$ that does not pass through any decorated vertices.      Then we have a diagonal embedding  $X^{{\rm Bound}_G/\sim} \to X^{{\rm Bound}_G}$  where ${\rm Bound}_G$ is the set of bound half-edges.    Let $F(v)$ be the set of free half edges adjacent to $v$.  Then there is a corresponding embedding  $$X^{{\rm Bound}_G/\sim} \to \prod_{v \in V_G \text{ decorated } }  X^{n(v)- F(v)}.$$

Let $$\prod_{v \in V_G} \cM'(v) \to \prod_{v \in V_G {\rm~ decorated}} X^{n(v)} \to \prod_{v \in V_G {\rm~ decorated}}  X^{n(v) - F(v)}$$  be the composite of the evaluation map with the map that forgets the points of $X$ corresponding to free edges, and let $\tilde \rS(G)'$ be the fiber product of this map with $X^{{\rm Bound}_G}.$

\begin{prop}\label{altdescription}
      If $G$ has a decorated vertex, then $\tilde \rS(G) \iso \tilde \rS(G)'$. Otherwise $\tilde \rS(G) \iso \tilde \rS(G)' \times X.$
\end{prop}
\begin{proof}
  By definition, a point of $\tilde S(G)$ consists of a pointed curve $C_v$ and a morphism  $f_v: C_v \to X$ for each $v \in G$,  such that $(C_v, f_v) \in \cM(v)$ and if $h_1, h_2 \in {\rm Half}(G)$ form an edge, then $f_{v_1}(p_1) = f_{v_2}(p_2)$, where $v_i,  p_i \in C_{v_i}$  are the vertices and marked points corresponding to $h_i$, for $i = 1,2$.  

If every vertex of $G$ is undecorated, then every $f_v$ is constant with the same value $x \in X$. So this data is equivalent to a collection of $C_v \in \cM'(v)$  and an element $x \in X$.   So we have that $\tilde \rS(G) \iso \tilde \rS(G)' \times X$.

If $G$ has a decorated vertex, then the value of $f_v$ for $v \in G$ undecorated is determined by the value of $f_w$ for $w \in G$ decorated: choose a path from $v$ to a decorated vertex that does not pass through any other decorated vertices.    Further, the choice of $f_w$ for $w$ decorated satisfies the condition that if $h_1$ and $h_2$ are non-free half-edges with $h_1 \sim h_2$   then $f_{v_1}(p_1) = f_{v_2}(p_2)$.  Conversely, a choice of maps $f_w \in \cM(w)$ for each $w$ undecorated satsifying the condition that if $h_1 \sim h_2$  then $f_{v_1}(p_1) = f_{v_2}(p_2)$ uniquely determines an element of $\tilde \rS(G)$.  So we have $\tilde \rS(G)' \iso \tilde \rS(G)$. 
\end{proof}

\begin{proof}[Proof of Proposition \ref{freebound}]
By Proposition \ref{altdescription} and K\"unneth, it suffices to prove that $H_i^{\BM}(\tilde S(G)') = 0$.  

Let  $v \in G$ be a vertex.    We define a morphism $\pi(v)$, with source $\cM'(v)$ casewise.  We will write $T(v)$ for the target of $\pi(v)$.  
\begin{enumerate}
\item  If $v$ is plain or has no free half edges,  let $\pi(v) = \id_{\cM'(v)}.$

\item Otherwise,  if $v$ is decorated and has at least one free half edge,  choose a free half edge $e$ of $v$ arbitrarily and  let $$\pi(v):\cM_{n(v), g(v)}(X,\alpha(v)) \to \cM_{1, g(v)}(X,\alpha(v))$$ be the morphism forgetting all the marked points except the one corresponding to $e$ 

\end{enumerate}

   In case $(2)$, $\pi(v)$ has fibers isomorphic to (a quotient by a finite group action of) $\Conf_{F(v)-1}(\Sigma_{g(v),1})$,  where $\Sigma_{g(v),1}$ is a once-punctured genus $g(v)$ surface.

Now consider the map $$\prod_{v} \pi(v): \prod_{v \in V_G} \cM'(v) \to \prod_{v \in V_G} T(v).$$  The map $$\prod_{v \in V_G} \cM'(v) \to \prod_{v ~ {\rm decorated }} X^{n(v) - F(v)}$$  factors through $\prod_{v} T(v)$.    So letting $Z$ be the fiber product of $\prod_{v \in G} T(v)$ with $$X^t \to \prod_{v \in V_G \text{ decorated } }  X^{n(v)- F(v)},$$   we obtain a map  $\pi: \tilde \rS(G)' \to Z$.

By the Leray spectral sequence for compactly supported cohomology and proper base change, it suffices to show that the fibers $\pi\inv(z), z \in Z$  satisfy $\rH_j^{\BM}(\pi\inv(z)) = 0$  for all $j$ such that $j  + g + (L \cdot \alpha) \leq \sum_{v \in V_G} F(v)$ .  The fibers of $\pi$ are homeomorphic to (a quotient by a finite group of) a product of configuration spaces:   a factor of $\Conf_{F(v) - 1}(\Sigma_{g(v),1})$ for each vertex satisfying case (2).  For $m \in \bbN$  we have that  $$\rH_i^{\BM}(\Conf_{m}(\Sigma_{g,1})) = 0$$ for all $i \leq m$ (see \cite[proof of Proposition 5.6]{tosteson2018stability}).  So by the K\"unneth formula  applied to the product of configuration spaces over vertices,  we see that $\rH_j^{\BM}(\pi\inv(z)) = 0$  for all $j \leq \sum_{v \in V_G} \max(F(v)-1, 0)$. Because the number of decorated vertices is $\leq L \cdot \alpha$ and the number of vertices with $g(v) >1 $ is $\leq g$,  we have that  $$\sum_{v \in V_G} \max(F(v)-1, 0) \geq \sum_{v \in V_G} F(v) - g - (L \cdot \alpha).$$  So if $j + g + (L \cdot \alpha) \leq \sum_{v \in V_G} F(v)$  then $\rH_j^{\BM}(\pi\inv(z)) = 0$ as desired.  
\end{proof}

\section{Bounds on Graphs}
This section contains two statements, which we will use to bound the number of graphs with certain combinatorial properties.  
\begin{prop}\label{maingraphbound}
	Let $\alpha \in \rN_1(X)$ be a sum of curve classes and $L \in \rN^1(X)$ be an ample class.  Let $g \in \bbN$  and $k \in \bbQ,  k \geq 0$.  Let $G$ be a stable decorated graph of genus $g$ and class $\alpha$.   Suppose that 
\begin{enumerate}
\item  the number of free half-edges of $G$ is $\leq i$,
\item no two distinct plain vertices of $G$ are adjacent.
\end{enumerate}
Then $G$ has $\leq \max(  i + 2g  + 2 (L \cdot \alpha), 1)$ plain vertices. 
\end{prop}
\begin{proof}
First we bound the number of bound half-edges of $G$.   To so, we assume that $G$ has at least one bound half-edge.
Note that there are at most $L \cdot \alpha$ decorated vertices of $G$, since every nonzero class $\beta$ which is a sum of curve classes must have $L \cdot \beta \geq 1$.   Now consider the following bipartite graph $G'$ associated to $G$: 
\begin{itemize} \item there is a vertex of $G'$ for each decorated vertex of $G$,
\item there is a vertex of $G'$  for each equivalence class of the equivalence relation $\sim$ on the set of bound half-edges of $G$,
\item there is an edge of $G'$ for each bound half-edge $h$ of $G$, connecting the vertices which correspond to the equivalence class of $\sim$ containing $h$ and to the vertex of $G$ that $h$ is adjacent to. 
\end{itemize}
Topologically, $G'$ is homeomorphic to a quotient of $G$,  so we have $\chi(G') \geq 1 - g$.  We may compute the Euler characteristic of $G'$ as  $$\sum_{v \in G \text{ decorated} } 1 + \sum_{h \in \rH(G) \text{ bound} } \left( \frac{1}{|[h]|} - 1 \right),$$ where $|[h]|$ denotes the size of the equivalence class of $\sim$  corresponding to $h$.   By the definition of bound edges, each equivalence class has size $\geq 2$.   So $$\chi(G') \leq |\{ v \in G \text{ decorated}\}| - 1/2 |\{ h \in \rH(G)  \text{ bound} \}|.$$
Using our bound on decorated vertices, and the inequality, $\chi(G') \geq 1-g$ we rearrange to obtain  $$|\{ h \in \rH(G)  \text{ bound} \}| \leq 2g - 2 + 2 (L \cdot \alpha).$$ 

Now we bound the number of plain vertices. Since every half-edge that is adjacent to a non-plain vertex is either free or bound we have $$\sum_{v \in G \text{ not plain}} n(v) \leq i + \max(2g - 2 + 2 (L \cdot \alpha), 0) \leq i + 2g + 2(L \cdot \alpha).$$
  Either $G$ consists of a single plain vertex, or every plain vertex is adjacent to a non-plain vertex.  In the second case, there is a half-edge associated to each plain vertex so there are $\leq i + 2g  + 2 (L \cdot \alpha)$  plain vertices.   So in total, there are $\leq\max(  i + 2g  + 2 (L \cdot \alpha),1)$ plain vertices.  
\end{proof}

The following proposition is similar to a specialization of  \cite[Lemma 4.2]{tosteson2018stability}, but with a slightly better bound. 
\begin{prop}\label{Moncombinatorics}
Let $G$ be a simply connected graph, all of whose vertices have valence $\geq 3$.  Let $i \in \bbN$, and suppose that $\sum_{v \in G} (n(v) - 3) = i$.  We further assume
\begin{enumerate}
\item  no trivalent vertex of $G$ has $\geq 2$ external edges,
 \item there are no adjacent trivalent vertices of $G$ such that both vertices have an external edge.
\end{enumerate}
Then the number of external edges of $G$ is $\leq 13i/2$.  
\end{prop}
\begin{proof}  
We assume that $G$ has more than one vertex, since if $G$ has a single vertex, its valence $n$ is necessarily $>3$, so that $i = n -3 > 0$, and hence  $n = i +3 \leq 13i/2$.  

For a vertex $v$ we write $\re\rx(v)$ for the number of external edges of $v$.  Then  we have that\begin{equation}
\sum_{v \in V_G} \left(1 - \frac{n(v) - \re\rx(v)}{2} \right) = \chi(G) = 1. \label{eulercharacteristic}
\end{equation}

  Rearranging we see that $$\sum_{v \in V_G} \re\rx(v) = 2 + |V_G| + \sum_{v \in V_G} (n(v) - 3).$$

To bound $|V_G|$,  we let $m_{3,1}$   (resp. $m_{3,0}$)  denote the number  of trivalent vertices of $G$ with one external edge  (resp. zero external edges).    Let $s$ be the number of vertices of valence $> 3$. 

We claim that $m_{3,0} \leq s -2$.  Indeed, in the  sum \eqref{eulercharacteristic}, each trivalent vertex with no external edges contributes $-1/2$.  Every other vertex contributes at most $1/2$, and because there are no trivalent vertices with $2$ external edges,  all of these vertices must be $>3$-valent and there must be at least $m_{3,0} + 2$ of them.  This establishes the inequality.  

Next, we claim that $2 m_{3,1} \leq  3 m_{3,0} + 3s + i$. Indeed, every internal half-edge adjacent to a trivalent vertex with one external edge must either (I) be adjacent to  a trivalent vertex with no external edges or (II) be adjacent to  vertex of valence $> 3$. There are at most $3 m_{3,0}$ edges satisfying case (I).    There are at most $$\sum_{v \in V_G, n(v) > 3} n(v) = 3s + \sum_{v \in V_G} (n(v) -3) \leq 3s +i$$  satisfying case (II).    This establishes the claim. 

Finally, because the vertices of $G$ are at least trivalent  the inequality $\sum_{v} n(v) - 3 \leq i$ implies that $s \leq i$.    Therefore we have that $$|V_G| = s + m_{3,0} + m_{3,1} \leq i + (i - 2) + \frac{3(i-2) + 3i + i}{2}.$$ Using this inequality, we get that $\sum_{v \in V_G} \re\rx(v) \leq 13i/2$.  
\end{proof}

\section{$\FSop$ module structures on genus $0$ curves}\label{genus0}

Let $\overline \cM_{e,n,0} \subseteq \overline \cM_{e,n}$  be the closed subvariety consisting of stable genus $e$ curves with $n$ marked points such that every irreducible component has genus $0$.   The subspace $\overline \cM_{e,n,0}$ is preserved by the action of $\BT\op$ and $H_i(\overline \cM_{e,n,0})$ is an $\FSop$ module.  
In this section, we will show that this $\FSop$ module  is finitely generated.    Our first proposition concerns the case $e = 0$.   

\begin{prop}\label{prop:genuszero}
	The $\FSop$ module $H_i(\overline \cM_{0,n})$ is generated in degrees $\leq \max(13 i/2 ,3) $ and has height $\leq \max(13 i/2,1)$.  
\end{prop}
\begin{proof}
	We treat the case $i = 0$ separately.   In this case $H_0(\overline \cM_{0,n}) = \bbQ$ for $n \geq 3$ and $0$ otherwise, with all transition maps the identity.  So the $\FSop$ module is generated in degree $3$ and a submodule of the free module generated in degree $1$.  

		The homology of $\overline \cM_{0,n}$ is additively generated by classes corresponding to undecorated, genus $0$ stable graphs $H$ with $n$ external edges and $i = \sum_{v \in H} (n(v) - 3)$.  This classes are subject to several relations, including the WDVV relation  (see e.g. the discussion in \cite[\S 5]{tosteson2018stability}).   Let $H \in \Stab(0,n)$  be a graph, and $[H] \in H_{i}(\bMon)$ be the corresponding class.  Suppose that $H$ has $> 13 i/2$ external edges.  Then by Propostion \ref{Moncombinatorics} we have that  $H$ either (1) has a trivalent vertex with $\geq 2$ external edges, or (2) contains $2$ adjacent trivalent vertices both of which have an external edge.  

In case (1), let $i,j \in [n]$ be the labels of the external edges of the trivalent vertex of $H$. Considering the surjection $\delta: [n] \to [n]/(i \sim j)$,  we see that $[H] = \delta^* [H_0]$, where $H_0$ is the graph obtained from $H$ by removing the trivalent vertex.

  In case (2), the WDVV relation implies that $[H] = [H']$ where $H'$ contains a trivalent vertex with two external edges.   So by the same argument, there is a $\delta
':[n] \to [n]/(i' \sim j')$  and an $H'_0$ with $\delta'{^*} [H_0'] = [H'] = [H]$.  

Thus if $n > 13i/2 $ every class of $H_i(\overline \cM_{0,n})$ is a linear combination of classes pushed forward from $H_i(\overline \cM_{0,n-1})$ under the action of $\FSop$.    So $H_i(\overline \cM_{0,n})$ is generated in degree $\leq 13i/2$.  
\end{proof}

\subsection{Resolution}   In order to show that $H_i(\overline \cM_{e,n,0})$  is a finitely generated $\FSop$ module,  we will express  the homology of  $\overline \cM_{e,n,0}$ in terms of the homology of $\overline \cM_{0,n}$ and gluing maps.  To do this we resolve $\overline \cM_{e,n,0}$, using the combinatorics of a category of stable graphs.   The motivation for this resolution is that there is a surjective gluing map $ \overline \cM_{0, n+2e} \to \overline \cM_{e,n,0}$  given by gluing together the $e$ pairs of marked points  $\{p_{2i-1},p_{2i}\}_{i = 1}^{e}$.   This map factors though $\overline \cM_{0, n + 2e}/( \rS_{e} \wr \bbZ/2)$, but the resulting quotient map is not an isomorphism,  because it is not one-to-one over the boundary strata.  To account for the discrepancy on boundary strata, we will define a functor from a category of stable graphs that maps the rose graph to $\overline\cM_{0, n + 2e}$. Instead of quotienting by $\rS_e \wr \bbZ/2$  (the automorphism group of the rose graph),  we will take the colimit of this functor to obtain $\overline \cM_{e,n,0}$.  


 Let $\cG_e$ be a skeleton of the category of connected graphs of genus $e$ with all vertices of genus $0$,  and contractions between them.   Thus $\cG_e$ contains a unique object for each isomorphism class of such graph.  (For our purposes, could take $\cG_e$ to be the category of all graphs,  but for notational simplicity and to avoid size issues,  we work with a skeleton instead).  

\begin{defn} For a graph $G \in \cG_e$,   we define   $$\cF_n(G) := \bigsqcup_{f:[n] \to \rV(G) \text{ stable }} \prod_{v \in \rV(G)}  \overline \cM_{0, f\inv(v) \sqcup n(v)},$$ where the disjoint union is over all functions $f: [n] \to \rV(G)$  such that $|f\inv(v)| + n(v) \geq 3$.    
\end{defn}

The space $\cF_n(G)$ parameterizes (isomorphism classes of) the following data:  
\begin{itemize} \item a nodal curve $C$ with $H^1(C, \cO_C) = 0$  
\item a collection of smooth points in $C$  marked by $[n]$  and  $\rH\ra\rl\rf(G)$
\end{itemize}
satisfying the conditions that
\begin{itemize}
\item  every connected component of $C$ is stable, has at least $3$ marked points, and at least one point marked by a half edge 
\item  two half edges are adjacent to the same vertex of $G$ if and only if their corresponding marked points lie on the same connected component of $C$.
\end{itemize}
Notice that the automorphism group of every such marked curve trivial.  

In terms of this description, we see that a contraction of graphs $G_1 \to G_2 \in \cG_e$  determines morphism $\cF_n(G_1) \to \cF_n(G_2)$  by gluing and relabelling points according to the contraction, as follows.   Let  $\iota: \rH\ra\rl\rf(G_2) \to \rH\ra\rl\rf(G_1)$ is the associated embedding of the uncontracted half edges.    The map $\cF_n(G_1) \to \cF_n(G_2)$ takes a curve $C \in \cF_n(G_1)$  to the curve  $$\frac{C}{p_h \sim p_{\sigma(h)}\text{  for all $h \not \in \im \iota$}},$$  and associates to $f \in \rH\ra\rl\rf(G_2)$  the labelled point $p_{\iota(f)}$ (considered as a point in the quotient curve).  

Thus $G \mapsto \cF_n(G)$ is functor from $\cG_e$ to varieties.  Further, for every $G \in \cG_e$ there is a gluing morphism ${\rm glue} : \cF_n(G) \to \overline \cM_{e,n,0}$, given by $${\rm glue}(C)  := C/\{p_h \sim p_{\sigma(h)}\}_{h \in \rH\ra\rl\rf(G)}.$$
These morphisms are compatible with $\cF_n(G_1) \to \cF_n(G_2)$, hence define a cone on the diagram $\cF_n$. We will show that the diagram $\cF_n$  determines the rational cohomology of the moduli space $\overline \cM_{e,n,0}$.  

\begin{remark} There are different ways of describing the relationship between $\cF_n$ and $\overline \cM_{e,n,0}$.   From a homotopical perspective, the homotopy type underlying (the stack corresponding to) $\overline \cM_{e,n,0}$  is the homotopy colimit of the diagram $\cF_n$.   From an algebro-geometric perspective,  the diagram of topoi  $(G \in \cG_e) \mapsto \Sh(\cF_n(G))$  augmented by $\Sh(\cF_n(G)) \to \Sh(\overline \cM_{e,n,0})$ is of cohomological descent.   These statements are not equivalent, but express similar ideas and have the same consequences for rational cohomology.   We will not prove either; instead we will opt for a more direct simplicial approach.  However it may be helpful for a reader with the relevant background to have these statements in mind.  
\end{remark}

Let $\rB(\cF_n)$  be the bar construction of $\cF_n$.   By definition $\rB(\cF_n)$ is a simplicial algebraic space, whose space of $r$ simplices is $$\rB(\cF_n)_r := \bigsqcup_{G_0 \to G_1 \to \dots \to G_{r}  \in \cG_e}  \cF_n(G_0).$$ The face and degeneracy maps of $\rB(\cF_n)$ are given in the usual way.    The maps $\cF_n(G) \to  \overline \cM_{e,n,0}$  yield an augmentation $\rB(\cF_n) \to \overline \cM_{e,n,0}$.  We define the homology of $\rB(\cF_n)$ to be the homology of the associated simplicial chain complex $r \mapsto \rC_*(\rB(\cF_n)_r)$.  

\begin{prop}
	The augmentation map  $\rB(\cF_n) \to \overline \cM_{e,n,0}$  induces an isomorphism on rational homology.  
\end{prop}
\begin{proof}
We first show that the augmentation map induces homology isomorphisms on fibers.  
	 Let $[C] \in \overline \cM_{e,n,0}$ be an isomorphism class of curve (a point in the coarse moduli space).  Then the fiber over $[C]$ is the simplicial algebraic space with $r$ simplices \begin{equation}r \mapsto \{  G_0 \to G_1 \to \dots \to G_r, D \in \cF_n(G_0), ~  [{\rm glue}(D)] =  [C]\}. \label{fibers} \end{equation}  This is the quotient by $\Aut(C)$ of \begin{equation} \label{chains1} r \mapsto \{  G_0 \to G_1 \to \dots \to G_r, D \in \cF_n(G_0), ~  \phi: \glue(D) \iso C\}.\end{equation}  Now, let $\tilde C$ be the normalization of $C$, and fix a graph $G_C \in \cG_e$ together with an isomorphism to the dual graph of $C$.  Under this isomorphism, $\tilde C$ is canonically an element of $\cF_n(G_C)$.   There is a canonical isomorphism $\glue(\tilde C) \iso C$.    We observe the following  facts.
\begin{enumerate}
\item  The simplicial set \eqref{chains1} is equivalent to the nerve of the category $\cG_e(C)$ whose objects are triples $(G \in \cG_{e}$,  $D \in \cF_n(G)$,   $\phi: \glue(D) \iso C)$  and whose morphisms are $G \to G' \in \cG_e$  such that $D \in \cF_n(G)$ maps to  $D' \in \cF_n(G')$  and $\phi'$ composed with the canonical map $\glue(D) \to \glue(D')$ equals $\phi$.  
\item   The object  $(G_C, \tilde C, \glue(\tilde C) \iso C)$  is initial in $\cG_e(C)$.  
\end{enumerate}
It follows that the simplicial set \eqref{chains1} deformation retracts to the vertex  $(G_C,  \tilde C \in \cF_n(G_C),  \glue (\tilde C) \iso C)$.  Quotienting by the action of $\Aut(C)$, it follows that augmentation map from the simplicial set \eqref{fibers} to a point induces an isomorphism on rational homology.  

The proposition follows from the statement on fibers, together with the Leray spectral sequence and the proper base change theorem.  Note that to show  $\rB(\cF_n) \to \overline \cM_{e,n,0}$ induces a rational homology isomorphism in degrees $\leq i$, it suffices to show the same for the restriction to the  $(i+1)$-skeleton $\rs\rk_{i+1}\rB(\cF_n) \to \overline \cM_{e,n,0}$.   The map $|\rs\rk_{i+1} \rB(\cF_n)| \to \overline \cM_{e,n,0}$  is proper, and the homology of the fibers vanishes in degree $\leq i$.  Hence by proper base change and the Leray spectral sequence, this map induces an isomorphism on homology in degree $\leq i$.  
\end{proof}

\subsection{Filtration} Our next aim is to construct a small chain complex computing the rational cohomology of $\rB(\cF_n)$.    To do so, we filter $\rB(\cF_n)$  using the canonical filtration on $\cG_e$ by number of edges:  let $F_i \rB(\cF_n)$ be the sub-simplicial space spanned by chains  $G_0 \to \dots \to G_r$  such that $G_0$ has $\leq e + i$ edges.   

We compute the relative homology of $(F_i \rB(\cF_n) , F_{i-1} \rB(\cF_n))$ as follows.    For a graph $G$ in $\cG_e$  let  $\rN(G/\cG_e)$ be the nerve of its overcategory: the simplicial set whose $r$ simplices are chains $G \onto G_0 \onto \dots \onto G_r$.  Let  $\rN(G/\cG_e^{> G})$ denote the subcomplex consisting of chains such that $G \to G_0$ contracts at least one edge.    Note that $\Aut(G)$  acts on both of these complexes.  

\begin{prop}\label{prop:associatedgraded}
		The relative rational homology $H_*(F_i \rB(\cF_n), F_{i-1} \rB(\cF_n))$  is naturally isomorphic to  $$\bigoplus_{G \in \cG_e, ~ |\rE\rd\rg\re(G)| = i+e} H_*( \cF_n(G)) \underset{\Aut(G)} \otimes  H_*(\rN(G/\cG_e) ,\rN(G/\cG_e^{> G}))$$ (as usual, all groups are taken with $\bbQ$ coefficients).  
\end{prop}
\begin{proof}
The chain complex computing the homology $F_i \rB(\cF_n) $  is $$\bigoplus_{G_0 \to \dots  \to G_r}  \rC_*(\cF_n(G_0))$$  where $G_0$ has $\leq i + e$ edges.   Taking the quotient by the subcomplex computing the homology of $F_{i-1} \rB(\cF_n)$, we get a complex $D$ which is a sum over chains where $G_0$ has exactly $i + e$ edges.   Similarly homology of the pair $(\rN(G/\cG_e), \rN(G/\cG_e^{> G}))$ is computed by a complex $\rC(G)$ which in degree $r$ is spanned by $G \iso G_0 \to G_1 \to \dots  \to G_r$.  Then $\Aut(G)$ acts freely on this complex,  and there is a canonical isomorphism between $$\bigoplus_{G \i \cG_e, ~ |E_G| = i+ e}  \rC_*(\cF(G))  \otimes_{\Aut(G)} \rC(G)$$ and the quotient complex $D$ given  by mapping  the element $G \iso^{\phi} G_0 \to \dots  \to G_r$, $c \in \rC(G)$   to  $G_0 \to \dots  \to G_r, \phi_* c \in \rC(G_0)$.
\end{proof}

Fix a graph $G \in \cG_e$.  There is a canonical equivalence of categories between $G/\cG_e$  and the poset of subgraphs $C \subseteq G$ which do not contain a cycle.  The equivalence takes $f: G \to G'$ to the subgraph of edges contracted by $f$.  Therefore  $G/\cG_e^{> G}$  is equivalent to the poset of non-empty independent subsets of the edges of $G$;  in other words, the poset of faces of the \emph{independence complex of $G$} (i.e. the independence complex of the matroid associated to $G$).   Since the nerve construction preserves equivalences and the nerve of a face poset of a simplicial complex is homeomorphic to its barycentric subsivision, we obtain the following.
\begin{prop}
 $\rN(G/\cG_e^{>G})$ is homotopy equivalent to the independence complex of $G$.
\end{prop}

Independence complexes of arbitrary matroids are pure and shellable of dimension equal to the rank of the matroid \cite[Theorem 7.3.3]{Bjorner}.  Since $G/\cG_e$ contains an initial object, $\rN(G/\cG_e)$ is contractible hence by the long exact sequence for a pair we obtain.  


\begin{prop}\label{relativecomputation}
	The homology of $(\rN(G/\cG_e) ,\rN(G/\cG_e^{> G}))$  is concentrated in degree $|{\rm Edge}(G)| - e$  and is of the same rank as the reduced homology of the independence complex of $G$.  
\end{prop}

Putting Propositions \ref{prop:associatedgraded} and \ref{relativecomputation}  together, we obtain the following spectral sequence associated to the filtration $F_i$.

\begin{prop}\label{hocolimSS}
	For a graph $G \in \cG_e$,  let $\rI(G)$ be the degree $E_G - e -1$ reduced (rational) homology of its independence complex.  There is a spectral sequence  $E^1_{p,q}(n) \implies H_{p+q}  (\overline \cM_{e,n,0}),$  with $$E^1_{p,q}(n) = \bigoplus_{G \in \cG_e, ~ |E_G|- e = p} H_q(\cF_n(G)) \underset{\Aut(G)} \otimes \rI(G).$$ As we vary $n \in \bbN$, these spectral sequences assemble to form a spectral sequence of $\FSop$ modules.  
\end{prop}  

Note that the final statement follows from the fact that $G \mapsto ([n] \mapsto \cF_n(G))$  defines a functor from $\cG_e \to \Var^{\FSop}$  and the filtration $F_i$ and the isomorphism of Proposition \ref{prop:associatedgraded}  are compatible with the action of $\FSop$.   Although we will not need this fact,  we note that Deligne's theory of weights implies that the spectral sequence of Proposition \ref{hocolimSS}  degenerates after the $E^1$ page.  From this spectral sequence, we obtain  the following.

\begin{prop}\label{genus0eloops}
		The $\FSop$ module  $n \mapsto H_i(\overline \cM_{e,n,0})$ is of height $\leq \max(13 i/2,1)$.  
\end{prop}
\begin{proof}
	 We prove that the $\FSop$ module $n \mapsto \bigoplus_{p+q = i} E^1_{p,q}(n)$ is of height $\leq \max(13i/2,1)$.    There are finitely many isomorphism classes of graphs with $\leq i + e$ edges, and $I(G)$ is finite dimensional for every graph $G$.  Therefore it suffices to show that for every graph $G \in \cG_e$ with $p+e$ edges  the $\FSop$ module  $n \mapsto H_{i - p}(\cF_n(G))$ is of height $\leq \max(13i/2,1)$.  Let $j = i - p$. By K\"unneth  we have  $$H_j(\cF_n(G)) = \bigoplus_{j_{v},~ v \in {\rm Vert}(G),~ \sum_v {j_v} = j}~ \bigoplus_{f: [n] \to {\rm Vert}(G)}~ \bigotimes_{v \in {\rm Vert}(G)}  H_{j_v} (\overline \cM_{0,f\inv(v) \sqcup n(v)}).$$   By Proposition \ref{shiftclosure} and Proposition \ref{prop:genuszero}, the module   $m \mapsto H_{j_v} (\overline \cM_{0, [m] \sqcup n(v)})$ is of height $\leq \max(13j_v/2,1) $.  Hence by Proposition \ref{convolutionclosure}  $n \mapsto H_j(\cF_n(G))$ is of height $$\leq \max_{j_v,~j = \sum_{v \in V_G} j_v}\left( \sum_{v \in G} \max(13j_v/2,1)\right)  =  \max(13 j/2 + |V_G| - 1, 1).$$   Since the number of vertices of $G$ is $p+1$, we have that the left hand side is $\max(13(i - p)/2 + p, 1) \leq \max(13i/2,1)$ as desired.  
\end{proof}

\begin{remark}
		The construction of the resolution $\rB(\cF_n) \to \overline \cM_{e,n,0}$  and the associated spectral sequence is general.  We record the essential properties of $\cG_e$ and the augmented diagram $\cF_n$ that we have used.
\begin{itemize}
\item  $\cG_e$ is a small EI-category, with finite automorphism groups, such that all ascending chains in the associated poset of isomorphism classes have finite length.  
\item The augmentation maps $\cF_n(G) \to \overline \cM_{e,n,0}$ are proper,  and for $C \in \overline \cM_{e,n,0}$ the category formed by the pairs $G \in \cG_e$  and $D \in\cF_n(G) \times_{ \overline \cM_{e,n,0}} C$   is $\bbQ$ acyclic.  
\end{itemize} 
\end{remark}

\section{The $\FSop$ action on $\overline \cM_{g,n}(X,\alpha)$}

Consider  the stratification of $\overline \cM_{g,n}(X, \alpha)$ by $\rQ(g,n, \alpha)$.  $\BT\op$ acts on $\overline \cM_{g,n}(X, \alpha)$ in a way that is compatible with the stratification,  inducing an action of $\FSop$ on the posets $\rQ(g,n, \alpha).$   (See \cite[Definition 3.14]{tosteson2018stability} for  a description of this action in the case $X =* $).

The $\FSop$ set $n \mapsto \rQ(g,n,\alpha)$ decomposes into orbits.  There is an orbit for every stable graph $K$ which is \emph{reduced}  in the sense of satisfying the condition:
\begin{itemize}
	\item  every plain vertex of $K$ has $\leq 1$ external edge.  
\end{itemize} 
 To see this,  note that we may associate to any $G \in \rQ(g,n,\alpha)$  a reduced graph $\tilde G$ in the same orbit, obtained by removing all but one external edge for each plain vertex of $G$, then  stabilizing and relabelling arbitrarily.   Two reduced graphs generate the same orbit if and only if they become isomorphic after a relabelling of their external edges.

Given a reduced stable decorated graph $K$ with $n(K) = n$ and  total class $\alpha$,  consider the $\BT\op$ space spanned by its $\FSop$ orbit   $$m \mapsto \bigsqcup_{G \in \rQ(g,m,\alpha),~ G = f^* K,~  f \in \FS(m, n(K))} \tilde S(G) $$

We may decompose this space as a product \begin{equation}Y(K) \times \bigsqcup_{f: m \to n(K)}  \prod_{h \in n(K)}~  \cN_K(h, f\inv(h)), \label{productdecomp} \end{equation} where $Y(K)$ and $\cN_K$ are defined as follows.   If the unique vertex $v$ adjacent to $i$ is not plain, we define $\cN_K(h,x)$ to be $\overline \cM_{0, x \sqcup *}$ for $|x| > 1$  and a point otherwise.   If $v$ is plain with $e$ self edges, then we define $$\cN_K(h, x):= \overline \cM_{ e,~ x\sqcup (n(v) - h), 0},$$ if $|x| + n(v) -1 > 2$ and a point otherwise.  Finally, $Y(K)$ is a subspace of the product of $\cM(w)$ over all vertices $w \in V_G$  which are  not plain or not  adjacent to an external edge:    $Y(K)$ is the subspace consisting of $(C_w, f_w : C_w \to X)$  such that if $p \in C_w$  and $p' \in C_{w'}$  correspond to connected half-edges of $K$, then $f_w(p) = f_{w'}(p').$

\begin{ex}
The following image depicts a reduced stable graph $K$ and the stable graphs in the orbit of $K$.   \vspace{-.1in}
\begin{center}
\includegraphics[scale=.8]{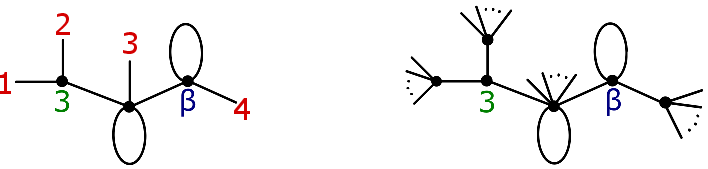} 
\vspace{-.1in}
\end{center}
Here $Y(K) = \cM_{3,3}(X,0) \times \cM_{0,4}(X,\beta)$  and $\cN_K(i,x) = \overline \cM_{0, x \sqcup  *}$ for $i = 1,2,4$  and $\cN_K(3,x) = \overline \cM_{1, x + 2, 0}$.  
\end{ex}

Using this decomposition, we prove the following proposition.

\begin{prop} \label{heightcomponent}
	Let $i \in \bbN$ and let $K$  be a reduced stable decorated graph.   Then the $\FSop$ module $$m \mapsto \bigoplus_{G \in \rQ(g,m,\alpha),~ G = f^* K,~  f \in \FS(m, n(K))} H_i^{\BM}(\tilde S(G))$$ has height $\leq \mainbound$.  
\end{prop}
\begin{proof}
The action of $\FSop$ on $\rQ(g,n,\alpha)$ preserves the number of free edges.    
	By Proposition \ref{freebound},  we may assume that $K$ has $\leq i + g + (\deg \alpha)$ free edges, since otherwise the $\FSop$ module is zero.  

By Proposition \ref{maingraphbound},  it follows that $K$ has $\leq \max(  i + 3g  + 3\deg \alpha), 1)$ plain vertices.  Since $K$ is reduced and there are $\leq i + g + \deg \alpha$ external edges adjacent to non-plain vertices, we have that $n(K) \leq  \max(2i + 4g  + 4  \deg \alpha , 1)$.    

Let $j \leq i$.  Then the $\FSop$ module $x \mapsto H_j^{\BM}(\cN(h,x))$ is of height $ \leq \max(13j/2,1) \leq (13i/2 + 1)$  by Proposition \ref{genus0eloops} and Proposition \ref{shiftclosure}.  From the decomposition \ref{productdecomp} the K\"unneth formula and Proposition \ref{convolutionclosure}, it follows that the $\FSop$ module of the Proposition has height  $\leq (13 i/2 + 1)(2i + 4g  + 4 \deg \alpha) + 1)$.  
\end{proof}  

Now we prove the main theorem.

\begin{proof}[Proof of Theorem \ref{mainthm}]
We apply the Borel--Moore homology spectral sequence associated to the stratification by $ \rQ(g,n,\alpha)$ and Proposition \ref{SSclosure}.   The group in degree $i$ is  $\bigoplus_{G \in \rQ(g,n, \alpha)}  H_i^{\BM}(S(G))$.  By Proposition \ref{quotienthomeo}, there is a surjection of $\FSop$ modules $$\bigoplus_{G \in \rQ(g,n,\alpha)} H_i^{\BM}(\tilde S(G)) \to \bigoplus_{G \in \rQ(g,n,\alpha)}  H_i^{\BM}(S(G)).$$

Decomposing according to the action of $\FSop$ on $\rQ(g,n,\alpha)$,  we have that the left hand side is $$\bigoplus_{[K] \text{ reduced  stable decorated} }  ~~~ \bigoplus_{G = f^* K,~  f \in \FS(n, n(K))} H_i^{\BM}(\tilde S(G)),$$ where the sum is over equivalence classes of reduced stable decorated graphs (with respect to relabelling the external edges).  By Proposition \ref{freebound}, we may restrict the sum to reduced stable graphs $K$ with $\leq i + g + \deg \alpha$ free edges,  since graphs with $>i$ free edges contribute $0$.    By Proposition \ref{heightcomponent},  each summand has height $\leq \mainbound$.  

So it suffices to show that there are finitely many terms in the sum.
\begin{lem}
	For $I \in \bbN$, there are finitely many reduced graphs in  $\bigsqcup_n \rQ(g,n,\alpha)$ with $\leq I$ free edges.  
\end{lem}
\begin{proof}
	By Proposition \ref{maingraphbound}  every stable decorated graph $K \in \sqcup_n \rQ(g,n,\alpha)$ with $\leq I$ free edges has $\leq \max(I + 2g  + 2 (L \cdot \alpha), 1)$ plain vertices.   Hence $K$ has $\leq \max(I + 2g  + 2 (L \cdot \alpha), 1) + i$  external edges.  Then it follows from the fact that $\rQ(g,n,\alpha) \subseteq \Stab(g,n,\alpha)$, which by Proposition \ref{prop:finite} is finite for every $n$.  
\end{proof} This completes the proof.  \end{proof}
\begin{proof}[Proof of Theorem \ref{thm:consequences}]
		The first part follows from \cite[Corollary 8.1.4]{sam2017grobner} and the second part follows from \cite[Theorem 1.17]{tosteson2021categorifications}.   
\end{proof}

\printbibliography

\end{document}